\input amstex

\input amstex
\input amsppt.sty
\magnification=\magstep1
\hsize=36truecc
 \vsize=23.5truecm
\baselineskip=14truept
 \NoBlackBoxes
\def\q{\quad}

\def\mod#1{\ (\text{\rm mod}\ #1)}

\def\t{\text}
\def\qtq#1{\q\t{#1}\q}
\def\mod#1{\ (\text{\rm mod}\ #1)}
\def\qtq#1{\q\t{#1}\q}
\def\f{\frac}
\def\e{\equiv}
\def\b{\binom}

\def\sls#1#2{(\f{#1}{#2})}
 \def\ls#1#2{\big(\f{#1}{#2}\big)}
\def\Ls#1#2{\Big(\f{#1}{#2}\Big)}
\par\q
\let \pro=\proclaim
\let \endpro=\endproclaim
\topmatter
\title Congruences involving $\binom{2k}k^2\binom{4k}{2k}m^{-k}$
\endtitle
\author ZHI-Hong Sun\endauthor
\affil School of the Mathematical Sciences, Huaiyin Normal
University,
\\ Huaian, Jiangsu 223001, PR China
\\ Email: zhihongsun$\@$yahoo.com
\\ Homepage: http://www.hytc.edu.cn/xsjl/szh
\endaffil

 \nologo \NoRunningHeads

\abstract{Let $p>3$ be a prime, and let $m$ be an integer with
$p\nmid m$. In the paper, by using the work of Ishii and Deuring's
theorem for elliptic curves with complex multiplication we solve
some conjectures of Zhi-Wei Sun concerning
$\sum_{k=0}^{p-1}\binom{2k}k^2\binom{4k}{2k}m^{-k}\mod {p^2}$.
  \par\q
}
 \endabstract
  \footnote"" {The author is
supported by the Natural Sciences Foundation of China (grant No.
10971078).}
\endtopmatter
\document
\subheading{1. Introduction}
\par For positive integers $a,b$ and $n$, if $n=ax^2+by^2$
for some integers $x$ and $y$, we briefly say that $n=ax^2+by^2$.
Let $p>3$ be a prime. In 2003, Rodriguez-Villegas[RV] posed some
conjectures on supercongruences modulo $p^2$. One of his conjectures
is equivalent to
$$\sum_{k=0}^{p-1}\f{\b{2k}k^2\b{4k}{2k}}{256^k}\e
\cases 4x^2-2p\mod{p^2}&\t{if $p=x^2+2y^2\e 1,3\mod 8$,}
\\0\mod{p^2}&\t{if $p\e 5,7\mod 8$.}
\endcases$$
  This conjecture has been solved  by
Mortenson[M] and Zhi-Wei Sun[Su2].
\par  Let $\Bbb Z$ be the set of
integers, and for a prime $p$ let $\Bbb Z_p$ be the set of rational
numbers whose denominator is coprime to $p$. Recently the author's
brother Zhi-Wei Sun[Su1] posed many conjectures for
$\sum_{k=0}^{p-1}\b{2k}k^2\b{4k}{2k}m^{-k}\mod{p^2}$,
 where $p>3$ is a prime and
$m\in\Bbb Z$ with $p\nmid m$. For example, he conjectured (see [Su1,
Conjecture A3])
$$\aligned\sum_{k=0}^{p-1}\f{\b{2k}k^2\b{4k}{2k}}{81^k}\e
\cases 0\mod {p^2}&\t{if $p\e 3,5,6\mod {7}$,}
\\4x^2-2p\mod {p^2}&\t{if $p\e 1,2,4\mod 7$ and so
$p=x^2+7y^2$.}
\endcases\endaligned\tag 1.1$$

\par Let $\{P_n(x)\}$ be the Legendre polynomials given by
 (see [MOS, pp.\;228-232], [G,
(3.132)-(3.133)])
$$P_n(x)=\f
1{2^n}\sum_{k=0}^{[n/2]}\b nk(-1)^k\b{2n-2k}nx^{n-2k} =\f 1{2^n\cdot
n!}\cdot\f{d^n}{dx^n}(x^2-1)^n,\tag 1.2$$ where $[a]$ is the
greatest integer not exceeding $a$.
 From (1.2)
we see that
$$P_n(-x)=(-1)^nP_n(x).\tag 1.3$$
\par Let $\sls am$ be the Jacobi symbol. For a prime $p>3$, In [S2] the author showed
that
$$ P_{[\f p4]}(t)\e
\sum_{k=0}^{[p/4]}\b{4k}{2k}\b{2k}k\Ls{1-t}{128}^k \e
-\Ls{6}p\sum_{x=0}^{p-1}\Ls{x^3-\f{3(3t+5)}2x+9t+7}p\mod
 p.\tag 1.4$$
 We
note that $p\mid\b{2k}k\b{4k}{2k}$ for $\f p4<k<p$.
 \par Let $p>3$
be a prime, $m\in\Bbb Z_p$, $m\not\e 0\mod p$ and
$t=\sqrt{1-256/m}$. In [S2] the author showed that
$$\sum_{k=0}^{p-1}\f{\b{2k}k^2\b{4k}{2k}}{m^k}\e
P_{[\f p4]}(t)^2\mod p.\tag 1.5$$ In the paper we show that
$$\sum_{k=0}^{p-1}\b{2k}k^2\b{4k}{2k}(x(1-64x))^k\e \Big(
\sum_{k=0}^{p-1}\b{2k}k\b{4k}{2k}x^k\Big)^2\mod {p^2}.\tag 1.6$$
 On
the basis of (1.5) and (1.6), we prove some congruences involving
$\sum_{k=0}^{p-1}\b{2k}k^2\b{4k}{2k}m^{-k}$ in the cases
$$m=81,648,12^4,28^4,1584^2,396^4, -144,-3969,
  -2^{10}\cdot 3^4,-2^{10}\cdot 21^4,
  -3\cdot 2^{12},-2^{14}\cdot 3^4\cdot 5.$$ Thus we partially solve
some conjectures posed by Zhi-Wei Sun in [Su1]. As examples, we
partially solve (1.1), and show that for primes $p\e \pm 1\mod 8$,
$$\sum_{k=0}^{p-1}\f{\b{2k}k^2\b{4k}{2k}}{1584^{2k}}
\e\cases 4x^2\mod p&\t{if $\sls p{11}=1$ and so $p=x^2+22y^2$,}
\\0\mod{p^2}&\t{if $\sls p{11}=-1$.}
\endcases$$

 \subheading{2. A general congruence modulo
$p^2$} \pro{Lemma 2.1} Let $m$ be a nonnegative integer. Then
$$\sum_{k=0}^m\b {2k}k^2\b{4k}{2k}\b k{m-k}(-64)^{m-k}
=\sum_{k=0}^m\b {2k}k\b{4k}{2k}\b{2(m-k)}{m-k}\b{4(m-k)}{2(m-k)}.$$
\endpro
\par We prove the lemma by using WZ method and Mathematica. Clearly
the result is true for $m=0,1$. Since both sides satisfy the same
recurrence relation $$\align 1024&(m+1)(2m+1)(2m+3)S(m)
    -8(2m+3)(8m^2+24m+19)S(m+1)\\&+(m+2)^3 S(m+2) = 0,\endalign$$
    we see that
    Lemma 2.1 is true.
 The  proof certificate for the left hand side is
$$ - \frac{4096 k^2(m+2)(m-2k)(m-2k+1)}{(m-k+1)(m-k+2)},$$
 and the proof certificate for the right hand side is
$$ \frac{16
k^2(4m-4k+1)(4m-4k+3)(16m^2-16mk+55m-26k+46)}{(m-k+1)^2(m-k+2)^2}.$$
\pro{Theorem 2.1} Let $p$ be an odd prime and let  $x$ be a
variable. Then
$$\sum_{k=0}^{p-1}\b{2k}k^2\b{4k}{2k}(x(1-64x))^k\e \Big(
\sum_{k=0}^{p-1}\b{2k}k\b{4k}{2k}x^k\Big)^2\mod {p^2}.$$
\endpro
Proof. It is clear that
$$\align &\sum_{k=0}^{p-1}\b{2k}k^2\b{4k}{2k}(x(1-64x))^k
\\&=\sum_{k=0}^{p-1}\b{2k}k^2\b{4k}{2k}x^k\sum_{r=0}^k\b kr(-64x)^r
\\&=\sum_{m=0}^{2(p-1)}x^m\sum_{k=0}^{min\{m,p-1\}}\b{2k}k^2\b{4k}{2k}\b
k{m-k}(-64)^{m-k}.\endalign$$
 Suppose $p\le m\le 2p-2$ and $0\le
k\le p-1$. If $k>\f p2$, then $p\mid \b{2k}k$ and so $p^2\mid
\b{2k}k^2$. If $k<\f p2$, then $m-k\ge p-k>k$ and so $\b k{m-k}=0$.
Thus, from the above and Lemma 2.1 we deduce
$$\align &\sum_{k=0}^{p-1}\b{2k}k^2\b{4k}{2k}(x(1-64x))^k
\\&\e \sum_{m=0}^{p-1}x^m\sum_{k=0}^m\b {2k}k^2\b{4k}{2k}
\b k{m-k}(-64)^{m-k}
\\&=\sum_{m=0}^{p-1}x^m\sum_{k=0}^m\b {2k}k\b{4k}{2k}
\b{2(m-k)}{m-k}\b{4(m-k)}{2(m-k)}
\\&=\sum_{k=0}^{p-1}\b{2k}k\b{4k}{2k}x^k\sum_{m=k}^{p-1}
\b{2(m-k)}{m-k}\b{4(m-k)}{2(m-k)}x^{m-k}
\\&=\sum_{k=0}^{p-1}\b{2k}k\b{4k}{2k}x^k
\sum_{r=0}^{p-1-k}\b{2r}r\b{4r}{2r}x^r
\\&=\sum_{k=0}^{p-1}\b{2k}k\b{4k}{2k}x^k\Big(\sum_{r=0}^{p-1}\b{2r}r
\b{4r}{2r}x^r -\sum_{r=p-k}^{p-1}\b{2r}r\b{4r}{2r}x^r\Big)
\\&=\Big(\sum_{k=0}^{p-1}\b{2k}k\b{4k}{2k}x^k\Big)^2
-\sum_{k=0}^{p-1}\b{2k}k\b{4k}{2k}x^k
\sum_{r=p-k}^{p-1}\b{2r}r\b{4r}{2r}x^r
 \mod{p^2}.\endalign$$
Now suppose $0\le k\le p-1$ and $p-k\le r\le p-1$. If $k\ge
\f{3p}4$, then $p^2\nmid (2k)!$, $p^3\mid (4k)!$ and so $
\b{2k}k\b{4k}{2k}=\f{(4k)!}{(2k)!k!^2}\e 0\mod{p^2}$. If $k < \f
p4$, then $r\ge p-k\ge \f{3p}4$ and so $
\b{2r}r\b{4r}{2r}=\f{(4r)!}{(2r)!r!^2}\e 0\mod{p^2}$. If $\f p4<k<\f
p2$, then $r\ge p-k>\f p2$, $p\nmid (2k)!$, $p\mid (4k)!$, $p\mid
\b{2r}r$ and $\b{2k}k\b{4k}{2k}=\f{(4k)!}{(2k)!k!^2}\e 0\mod p$. If
$\f p2<k<\f{3p}4$, then $r\ge p-k>\f p4$, $p\mid \b{2k}k$ and
$\b{2r}r\b{4r}{2r}=\f{(4r)!}{(2r)!r!^2}\e 0\mod p$. Hence we always
have $\b{2k}k\b{4k}{2k}\b{2r}r\b{4r}{2r}\e 0\mod{p^2}$ and so
$$\sum_{k=0}^{p-1}\b{2k}k\b{4k}{2k}x^k
\sum_{r=p-k}^{p-1}\b{2r}r\b{4r}{2r}x^r
 \e 0\mod{p^2}.$$
Now combining all the above we obtain the result.

 \pro{Corollary 2.1} Let $p>3$ be a
prime and $m\in\Bbb Z_p$ with $m\not\e 0\mod p$. Then
$$\sum_{k=0}^{p-1}\f{\b{2k}k^2\b{4k}{2k}}{m^k}
\e\Big(\sum_{k=0}^{p-1}\b{2k}k\b{4k}{2k}
\Big(\f{1-\sqrt{1-256/m}}{128}\Big)^k\Big)^2 \mod{p^2}.$$
\endpro
Proof. Taking $x=\f{1-\sqrt{1-256/m}}{128}$ in Theorem 2.1 we deduce
the result.
 \pro{Corollary 2.2} Let $p>3$ be a prime and $m\in\Bbb
Z_p$ with $m\not\e 0,256\mod p$. Then
$$\sum_{k=0}^{[p/4]}\f{\b{2k}k^2\b{4k}{2k}}{m^k}
\e 0\mod p\qtq{implies}\sum_{k=0}^{p-1}\f{\b{2k}k^2\b{4k}{2k}}{m^k}
\e 0\mod {p^2}.$$
\endpro
Proof. For $\f p4<k<p$ we see that
$\b{2k}k^2\b{4k}{2k}=\f{(4k)!}{k!^4}\e 0\mod p$.  Suppose
$\sum_{k=0}^{[p/4]}\f{\b{2k}k^2\b{4k}{2k}}{m^k} \e 0\mod p$. Then
$$\sum_{k=0}^{p-1}\f{\b{2k}k^2\b{4k}{2k}}{m^k} \e
\sum_{k=0}^{[p/4]}\f{\b{2k}k^2\b{4k}{2k}}{m^k}\e 0\mod p.$$ Using
Corollary 2.1 we see that
$$\sum_{k=0}^{p-1}\b{2k}k\b{4k}{2k}\Big(\f{1-\sqrt{1-256/m}}{128}\Big)^k
\e 0\mod p.$$ Thus the result follows from Corollary 2.1.

\pro{Corollary 2.3} Let $p\e 1,3\mod 8$ be a prime and $p=c^2+2d^2$
with $c,d\in\Bbb Z$ and $c\e 1\mod 4$. Then
$$\sum_{k=0}^{p-1}\f{\b{2k}k\b{4k}{2k}}{128^k}
\e (-1)^{[\f p8]+\f{p-1}2}\Big(2c-\f p{2c}\Big)\mod{p^2}.$$
\endpro
Proof. By [S2, Theorem 2.1] we have
$$\sum_{k=0}^{p-1}\f{\b{2k}k\b{4k}{2k}}{128^k}\e
\sum_{k=0}^{[p/4]}\f{\b{2k}k\b{4k}{2k}}{128^k}
 \e (-1)^{[\f p8]+\f{p-1}2}2c\mod p.$$
  Set $\sum_{k=0}^{p-1}\f{\b{2k}k\b{4k}{2k}}{128^k}
=(-1)^{[\f p8]+\f{p-1}2}2c+qp$. Then
 $$\Big(\sum_{k=0}^{p-1}\f{\b{2k}k\b{4k}{2k}}{128^k}\Big)^2
 =((-1)^{[\f p8]+\f{p-1}2}2c+qp)^2\e
 4c^2+(-1)^{[\f p8]+\f{p-1}2}4cqp\mod{p^2}.$$ Taking $x=\f 1{128}$
  in Theorem 2.1 we get
 $$\sum_{k=0}^{p-1}\f{\b{2k}k^2\b{4k}{2k}}{256^k}
 \e \Big(\sum_{k=0}^{p-1}\f{\b{2k}k\b{4k}{2k}}{128^k}\Big)^2\mod{p^2}.$$
 From [M] and [Su2] we have
 $$\sum_{k=0}^{p-1}\f{\b{2k}k^2\b{4k}{2k}}{256^k}\e 4c^2-2p\mod{p^2}.$$
 Thus
 $$4c^2-2p\e  \Big(\sum_{k=0}^{p-1}\f{\b{2k}k\b{4k}{2k}}{128^k}\Big)^2
 \e 4c^2+(-1)^{[\f p8]+\f{p-1}2}4cqp\mod{p^2}$$ and hence
 $q\e -(-1)^{[\f p8]+\f{p-1}2}\f 1{2c}\mod p$. So the corollary is proved.

 \subheading{3. Congruences for $\sum_{k=0}^{p-1}
 \b{2k}k^2\b{4k}{2k}m^{-k}$}
 \par Let $p>3$ be a prime and $m\in\Bbb Z$ with $p\nmid m$.
  In the section we partially solve Z.W.
Sun's conjectures on $\sum_{k=0}^{p-1}\b{2k}k^2 \b{4k}{2k}m^{-k}\mod
{p^2}$. \pro{Lemma 3.1 ([S2, (4.2)])} Let $p>3$ be a prime and let
$u$ be a variable. Then
$$P_{[\f p4]}(u)\e -\Ls 6p\sum_{x=0}^{p-1}\Big(x^3-\f
32(3u+5)x+9u+7\Big)^{\f{p-1}2}\mod p.$$
\endpro

 \pro{Lemma 3.2 ([S2, Theorem 4.1])} Let $p$ be
an odd prime, $m\in\Bbb Z_p$, $m\not\e 0\mod p$ and
$t=\sqrt{1-256/m}$. Then
$$\sum_{k=0}^{p-1}\f{\b{2k}k^2\b{4k}{2k}}{m^k}\e
P_{[\f p4]}(t)^2\e
\Big(\sum_{x=0}^{p-1}(x^3+4x^2+(2-2t)x)^{\f{p-1}2}\Big)^2\mod p.$$
\endpro
\pro{Lemma 3.3} Let $p$ be an odd prime, $m\in\Bbb Z_p$, $m\not\e
0\mod p$ and $t=\sqrt{1-256/m}$. If $P_{[\f p4}(t)\e 0\mod p$, then
$$\sum_{k=0}^{p-1}\f{\b{2k}k^2\b{4k}{2k}}{m^k}\e 0\mod{p^2}.$$
\endpro
Proof. By Corollary 2.1 we have
$$\sum_{k=0}^{p-1}\f{\b{2k}k^2\b{4k}{2k}}{m^k}\e
\Big(\sum_{k=0}^{p-1}\b{2k}k\b{4k}{2k}\Ls{1-t}{128}^k\Big)^2\mod{p^2}.
$$ Observe that $p\mid \b{2k}k\b{4k}{2k}$ for $\f p4<k<p$. From
[S2, Lemma 2.2] we see that
$$\sum_{k=0}^{p-1}\b{2k}k\b{3k}k\Ls{1-t}{128}^k
\e P_{[\f p4]}(t)\mod p.$$ Thus the result follows.

\pro{Lemma 3.4 ([S3, Lemma 4.1])} Let $p$ be an odd prime and let
$a,m,n$ be p-adic integers. Then
$$\sum_{x=0}^{p-1}(x^3+a^2mx+a^3n)^{\f{p-1}2}\e a^{\f{p-1}2}
\sum_{x=0}^{p-1}(x^3+mx+n)^{\f{p-1}2}\mod p.$$ Moreover, if $a,m,n$
are congruent to some integers, then
$$\sum_{x=0}^{p-1}\Ls{x^3+a^2mx+a^3n}p=
\Ls ap\sum_{x=0}^{p-1}\Ls{x^3+mx+n}p.$$
\endpro

\pro{Theorem 3.1} Let $p\not=2,3,7$ be a prime. Then
$$\aligned &P_{[\f p4]}\Ls{5\sqrt{-7}}9\\&\e\cases
-\sls {3(7+\sqrt{-7})}p\sls C72C\mod p&\t{if $p\e 1,2,4\mod 7$ and
so $p=C^2+7D^2$,}
\\0\mod p&\t{if $p\e 3,5,6\mod 7$}
\endcases\endaligned$$
and
$$\sum_{k=0}^{p-1}\f{\b{2k}k^2\b{4k}{2k}}{81^k}
\e\cases 4C^2\mod p
&\t{if $p\e 1,2,4\mod 7$ and so $p=C^2+7D^2$,}
\\0\mod {p^2}&\t{if $p\e 3,5,6\mod 7$.}
\endcases$$
\endpro
Proof. By Lemma 3.1 we have
$$P_{[\f p4]}\Ls{5\sqrt{-7}}9
\e -\Ls 6p\sum_{x=0}^{p-1}\Big(x^3-\f
52(3+\sqrt{-7})x+7+5\sqrt{-7}\Big)^{\f{p-1}2}\mod p.$$ Since
$$\f{-\f
52(3+\sqrt{-7})}{-35}=\Ls{1-\sqrt{-7}}{2\sqrt{-7}}^2\qtq{and}
\f{7+5\sqrt{-7}}{-98}=\Ls{1-\sqrt{-7}}{2\sqrt{-7}}^3,$$ by the above
and Lemma 3.4 we have
$$P_{[\f p4]}\Ls{5\sqrt{-7}}9
\e -\Ls
6p\Ls{1-\sqrt{-7}}{2\sqrt{-7}}^{\f{p-1}2}\sum_{x=0}^{p-1}\Ls{x^3-35x-98}p
\mod p.$$ Observe that $(x+7)^3-35(x+7)-98=x^3+21x^2+112x$. By the
work of Rajwade ([R1,R2]), we get
$$\aligned \sum_{x=0}^{p-1}\Ls{x^3-35x-98}p&
=\sum_{x=0}^{p-1}\Ls{x^3+21x^2+112x}p
\\&=\cases 2C\sls C7&\t{if $p=C^2+7D^2\e 1,2,4\mod 7$,}
\\0&\t{if $p\e 3,5,6\mod 7$.}
\endcases\endaligned\tag 3.1$$
For $p\e 1,2,4\mod 7$ we see that
$$\Ls 6p\Ls{1-\sqrt{-7}}{2\sqrt{-7}}^{\f{p-1}2}
=\Ls 6p\Ls{7+\sqrt{-7}}{2\cdot (-7)}^{\f{p-1}2}\e \Ls
3p\Ls{7+\sqrt{-7}}p\mod p.$$ Thus, from the above we deduce the
congruence for $P_{[\f p4]}\sls{5\sqrt{-7}}9\mod p$. Applying Lemmas
3.2 and 3.3 we obtain the remaining result.
\par\q
\par Let $p>3$ be a prime and let $\Bbb F_p$ be the field of $p$ elements.
For $m,n\in \Bbb F_p$ let $\#E_p(x^3+mx+n)$ be the number of points
on the curve $E_p$: $y^2=x^3+mx+n$ over the field $\Bbb F_p$. It is
well known that
$$\#E_p(x^3+mx+n)=p+1+\sum_{x=0}^{p-1}\Ls{x^3+mx+n}p.\tag 3.2$$
Let $K=\Bbb Q(\sqrt{-d})$ be an imaginary quadratic field and the
curve $y^2=x^3+mx+n$ has complex multiplication by $K$. By Deuring's
theorem ([C, Theorem 14.16],[PV],[I]), we have
$$\#E_p(x^3+mx+n)=\cases p+1&\t{if $p$ is inert in $K$,}
\\p+1-\pi-\bar{\pi}&\t{if $p=\pi\bar{\pi}$ in $K$,}
\endcases\tag 3.3$$ where $\pi$ is in an order in $K$ and
$\bar{\pi}$ is the conjugate number of $\pi$. If $4p=u^2+dv^2$ with
$u,v\in\Bbb Z$, we may take $\pi=\f 12(u+v\sqrt{-d})$. Thus,
$$\sum_{x=0}^{p-1}\Ls{x^3+mx+n}p=\cases \pm u&\t{if $4p=u^2+dv^2$
with $u,v\in\Bbb Z$,}\\0&\t{otherwise.}
\endcases\tag 3.4$$
In [JM] and [PV] the sign of $u$ in (3.4) was determined for those
imaginary quadratic fields $K$ with class number $1$. In [LM] and
[I] the sign of $u$ in (3.4) was determined for imaginary quadratic
fields $K$ with class number $2$.

\pro{Theorem 3.2} Let $p$ be a prime such that $p\e \pm 1\mod{12}$.
Then
$$P_{[\f p4]}\Big(\f 7{12}\sqrt 3\Big)\e
\cases\sls{2+2\sqrt 3}p2x
\mod p&\t{if $p=x^2+9y^2\e 1\mod{12}$ with $3\mid x-1$,} \\
0\mod p&\t{if $p\e 11\mod{12}$}\endcases
$$ and
$$\sum_{k=0}^{p-1}\f{\b{2k}k^2\b{4k}{2k}}{(-12288)^k}
\e\cases 4x^2\mod p&\t{if $p=x^2+9y^2\e 1\mod{12}$,}
\\0\mod{p^2}&\t{if $p\e 11\mod{12}$.}
\endcases$$
\endpro
Proof. From [I, p.133] we know that the elliptic curve defined by
the equation $y^2=x^3-(120+42\sqrt 3)x+448+336\sqrt 3$ has complex
multiplication by the order of discriminant $-36$. Thus, by (3.4)
and [I, Theorem 3.1] we have
$$\aligned&\sum_{n=0}^{p-1}\Ls{n^3-(120+42\sqrt 3)x+448+336\sqrt 3}
p\\&=\cases -2x\sls{1+\sqrt 3}p&\t{if $p=x^2+9y^2\e 1\mod{12}$ with
$3\mid x-1$,}
\\0&\t{if $p\e 11\mod{12}.$}\endcases\endaligned$$
By Lemma 3.1 we have
$$\align P_{[\f p4]}\Big(\f 7{12}\sqrt 3\Big)&\e
-\Ls 6p\sum_{n=0}^{p-1}\Big(n^3-\f{60+21\sqrt 3)}8n+\f{28+21\sqrt
3}4\Big)^{\f{p-1}2}
\\&\e -\Ls 6p\sum_{n=0}^{p-1}\Big(\ls n4^3-\f{60+21\sqrt 3)}8\cdot\f
n4+\f{28+21\sqrt 3}4\Big)^{\f{p-1}2}
\\&\e -\Ls 6p\sum_{n=0}^{p-1}\Ls{n^3-(120+42\sqrt 3)x+448+336\sqrt 3}
p\mod p.\endalign$$ Now combining all the above we obtain the
congruence for $P_{[\f p4]}(\f 7{12}\sqrt 3)\mod p$. Applying Lemmas
3.2 and 3.3 we deduce the remaining result.
\par\q
\newline{\bf Remark 3.1} In [Su1, Conjecture A24], Z.W. Sun
conjectured that for any prime $p>3$,
$$\align&\sum_{k=0}^{p-1}\f{\b{2k}k^2\b{4k}{2k}}{(-12288)^k}
\\&\e\cases (-1)^{[\f x6]}(4x^2-2p)\mod p&\t{if $p=x^2+y^2\e 1\mod{12}$
and $4\mid x-1$,}
\\-4\sls{xy}3xy\mod{p^2}&\t{if $p=x^2+y^2\e 5\mod{12}$ and $4\mid
x-1$,}
\\0\mod{p^2}&\t{if $p\e 3\mod{4}$.}
\endcases\endalign$$
\par If $p$ is a prime such that $p=x^2+5y^2\e 1,9\mod{20}$,
 by using [LM, Theorem 11] the author proved in [S2,
Theorem 4.7] that
$$\sum_{k=0}^{p-1}\f{\b{2k}k^2\b{4k}{2k}}{(-1024)^k}\e 4x^2\mod p.$$
Now we give similar results concerning $x^2+13y^2$ and $x^2+37y^2$.

 \pro{Theorem 3.3}
Let $p$ be an odd prime such that $p\not=3$ and $\sls {13}p=1$. Then
$$\sum_{k=0}^{p-1}\f{\b{2k}k^2\b{4k}{2k}}{(-82944)^k}
\e\cases 4x^2\mod p&\t{if $p=x^2+13y^2\e 1\mod 4$,}
\\0\mod{p^2}&\t{if $p\e 3\mod 4$.}
\endcases$$
\endpro
Proof. From [LM, Table II] we know that the elliptic curve defined
by the equation $y^2=x^3+4x^2+(2-\f 59\sqrt{13})x$ has complex
multiplication by the order of discriminant $-52$. Thus, by (3.4) we
have
$$\aligned&\sum_{n=0}^{p-1}\Ls{n^3+4n^2+(2-\f 59\sqrt{13})n}p
\\&=\cases 2x&\t{if $p\e 1\mod 4$ and so $p=x^2+13y^2$,}
\\0&\t{if $p\e 3\mod 4.$}\endcases\endaligned$$
Now taking $m=-2^{10}\cdot 3^4$ and $t=\f 5{18}\sqrt{13}$ in Lemmas
3.2 and 3.3 and applying the above we deduce the result.
\par\q
\newline{\bf Remark 3.2} Let $p\not=3,13$ be an odd prime. In [Su1, Conjecture A17],
Z.W. Sun conjectured that
$$\sum_{k=0}^{p-1}\f{\b{2k}k^2\b{4k}{2k}}{(-82944)^k}
\e\cases 4x^2-2p\mod {p^2}&\t{if $\sls {13}p=\sls{-1}p=1$ and so
$p=x^2+13y^2$,}\\2p-2x^2\mod {p^2}&\t{if $\sls {13}p=\sls{-1}p=-1$
and so $2p=x^2+13y^2$,}
\\0\mod{p^2}&\t{if $\sls {13}p=-\sls{-1}p$.}
\endcases$$

\pro{Theorem 3.4} Let $p$ be an odd prime such that $p\not=3,7$ and
$\sls {37}p=1$. Then
$$\sum_{k=0}^{p-1}\f{\b{2k}k^2\b{4k}{2k}}{(-2^{10}\cdot 21^4)^k}
\e\cases 4x^2\mod p&\t{if $p\e 1\mod 4$ and so $p=x^2+37y^2$,}
\\0\mod{p^2}&\t{if $p\e 3\mod 4$.}
\endcases$$
\endpro
Proof. From [LM, Table II] we know that the elliptic curve defined
by the equation $y^2=x^3+4x^2+(2-\f {145}{441}\sqrt{37})x$ has
complex multiplication by the order of discriminant $-148$. Thus, by
(3.4) we have
$$\aligned&\sum_{n=0}^{p-1}\Ls{n^3+4n^2+(2-\f {145}{441}\sqrt{37})n}p
\\&=\cases 2x&\t{if $p\e 1\mod 4$ and so $p=x^2+37y^2$,}
\\0&\t{if $p\e 3\mod 4.$}\endcases\endaligned$$
Now taking $m=-2^{10}\cdot 21^4$ and $t=\f {145}{882}\sqrt{37}$ in
Lemmas 3.2 and 3.3 and applying the above we deduce the result.
\par\q
\newline{\bf Remark 3.3} Let $p\not=3,7,37$ be a prime. In
 [Su1, Conjecture A19],
Z.W. Sun conjectured that
$$\sum_{k=0}^{p-1}\f{\b{2k}k^2\b{4k}{2k}}{(-2^{10}\cdot 21^4)^k}
\e\cases 4x^2-2p\mod {p^2}&\t{if $\sls {37}p=\sls{-1}p=1$ and so
$p=x^2+37y^2$,}\\2p-2x^2\mod {p^2}&\t{if $\sls {37}p=\sls{-1}p=-1$
and so $2p=x^2+37y^2$,}
\\0\mod{p^2}&\t{if $\sls {37}p=-\sls{-1}p$.}
\endcases$$

\par Let $b\in\{3,5,11,29\}$ and $f(b)=48^2,12^4,1584^2,396^4$ according as
$b=3,5,11,29$. For any odd prime $p$ with $p\nmid bf(b)$, Z.W. Sun
conjectured that ([Su1, Conjectures A14, A16, A18 and A21])
$$\sum_{k=0}^{p-1}\f{\b{2k}k^2\b{4k}{2k}}{f(b)^k}
\e\cases 4x^2-2p\mod {p^2}&\t{if $\sls 2p=\sls{-b}p=1$ and so
$p=x^2+2by^2$,}\\2p-8x^2\mod {p^2}&\t{if $\sls 2p=\sls{-b}p=-1$ and
so $p=2x^2+by^2$,}
\\0\mod{p^2}&\t{if $\sls 2p=-\sls{-b}p$.}
\endcases\tag 3.5$$
\par Now we partially solve the above conjecture.
\pro{Theorem 3.5} Let $p$ be an odd prime such that $p\e \pm 1 \mod
8$. Then
$$P_{[\f p4]}\Ls{2\sqrt 2}3\e
\cases (-1)^{\f{p-1}2}\sls{\sqrt 2}p\sls x32x\mod p&\t{if
$p=x^2+6y^2\e 1,7\mod {24}$,}
\\0\mod p&\t{if $p\e 17,23\mod {24}$}\endcases $$ and
$$\sum_{k=0}^{p-1}\f{\b{2k}k^2\b{4k}{2k}}{48^{2k}}
\e\cases 4x^2\mod p&\t{if $p=x^2+6y^2\e 1,7\mod {24}$,}
\\0\mod{p^2}&\t{if $p\e 17,23\mod {24}$.}
\endcases$$
\endpro
Proof. From [I, p.133] we know that the elliptic curve defined by
the equation $y^2=x^3+(-21+12\sqrt 2)x-28+22\sqrt 2$ has complex
multiplication by the order of discriminant $-24$. Thus, by (3.4)
and [I, Theorem 3.1] we have
$$\aligned&\sum_{n=0}^{p-1}\Ls{n^3+(-21+12\sqrt 2)n-28+22\sqrt 2}p
\\&=\cases 2x\sls{2x}3\sls{1+\sqrt 2}p
&\t{if $p\e 1,7\mod {24}$ and so $p=x^2+6y^2$,}
\\0&\t{if $p\e 17,23\mod {24}.$}\endcases\endaligned$$
By Lemma 3.1 we have
$$\align P_{[\f p4]}\Ls{2\sqrt 2}3&\e
-\Ls 6p\sum_{n=0}^{p-1}\Big(n^3-\f {15+6\sqrt 2}2n+7+6\sqrt
2\Big)^{\f{p-1}2}\mod p.\endalign$$ Since
$$\f{-(15+6\sqrt 2)/2}
{-21+12\sqrt 2}=\Ls{\sqrt 2+1}{\sqrt 2}^2\qtq{and} \f{7+6\sqrt
2}{-28+22\sqrt 2}=\Ls{\sqrt 2+1}{\sqrt 2}^3,$$ by Lemma 3.4 and the
above we have
$$\align P_{[\f p4]}\Ls{2\sqrt 2}3&\e
-\Ls 6p  \Ls{\sqrt 2(\sqrt
2+1)}p\sum_{n=0}^{p-1}\Ls{n^3+(-21+12\sqrt 2)n-28+22\sqrt 2}p
\\&\e\cases -\sls 6p\sls{\sqrt 2}p2x\sls{2x}3
\mod p&\t{if $p=x^2+6y^2\e 1,7\mod{24}$,}
\\0\mod p&\t{if $p\e 17,23\mod p$.}
\endcases\endalign$$
This yields the result for $P_{[\f p4]}\sls{2\sqrt 2}3\mod p$.
Taking $m=48^2$ and $t=\f 23\sqrt 2$ in Lemmas 3.2 and 3.3 and
applying the above we deduce the remaining result.

\pro{Theorem 3.6} Let $p$ be a prime such that $p\e \pm 1\mod 5$.
Then
$$\sum_{k=0}^{p-1}\f{\b{2k}k^2\b{4k}{2k}}{12^{4k}}
\e\cases 4x^2\mod p&\t{if $p=x^2+10y^2\e 1,9,11,19\mod {40}$,}
\\0\mod{p^2}&\t{if $p\e 21,29,31,39\mod {40}$.}
\endcases$$
\endpro
Proof. From [LM, Table II] we know that the elliptic curve defined
by the equation $y^2=x^3+4x^2+(2-\f 89\sqrt{5})x$ has complex
multiplication by the order of discriminant $-40$. Thus, by (3.4) we
have
$$\aligned&\sum_{n=0}^{p-1}\Ls{n^3+4n^2+(2-\f 89\sqrt 5)n}p
\\&=\cases 2x&\t{if $p\e 1,9,11,19\mod {40}$ and so $p=x^2+10y^2$,}
\\0&\t{if $p\e 21,29,31,39\mod {40}.$}\endcases\endaligned$$
Now taking $m=12^4$ and $t=\f 49\sqrt 5$ in Lemmas 3.2 and 3.3 and
applying the above we deduce the result.
\par\q

\pro{Theorem 3.7} Let $p$ be a prime such that $p\e \pm 1\mod 8$.
Then
$$\sum_{k=0}^{p-1}\f{\b{2k}k^2\b{4k}{2k}}{1584^{2k}}
\e\cases 4x^2\mod p&\t{if $\sls p{11}=1$ and so $p=x^2+22y^2$,}
\\0\mod{p^2}&\t{if $\sls p{11}=-1$.}
\endcases$$
\endpro
Proof. From [LM, Table II] we know that the elliptic curve defined
by the equation $y^2=x^3+4x^2+(2-\f {140}{99}\sqrt{2})x$ has complex
multiplication by the order of discriminant $-88$. Thus, by (3.4) we
have
$$\aligned&\sum_{n=0}^{p-1}\Ls{n^3+4n^2+(2-\f {140}{99}\sqrt{2})n}p
\\&=\cases 2x&\t{if $\sls p{11}=1$ and so $p=x^2+22y^2$,}
\\0&\t{if $\sls p{11}=-1.$}\endcases\endaligned$$
Now taking $m=1584^2$ and $t=\f {70}{99}\sqrt 2$ in Lemmas 3.2 and
3.3 and applying the above we deduce the result.
 \pro{Theorem 3.8}
Let $p$ be an odd prime such that $\sls {29}p=1$. Then
$$\sum_{k=0}^{p-1}\f{\b{2k}k^2\b{4k}{2k}}{396^{4k}}
\e\cases 4x^2\mod p&\t{if $p=x^2+58y^2\e 1,3\mod 8$,}
\\0\mod{p^2}&\t{if $p\e 5,7\mod 8$.}
\endcases$$
\endpro
Proof. From [LM, Table II] we know that the elliptic curve defined
by the equation $y^2=x^3+4x^2+(2-\f {3640}{9801}\sqrt{29})x$ has
complex multiplication by the order of discriminant $-232$. Thus, by
(3.4) we have
$$\aligned&\sum_{n=0}^{p-1}\Ls{n^3+4n^2+(2-\f {3640}{9801}\sqrt{29})n}p
\\&=\cases 2x&\t{if $p\e 1,3\mod 8$ and so $p=x^2+58y^2$,}
\\0&\t{if $p\e 5,7\mod 8.$}\endcases\endaligned$$
Now taking $m=396^4$ and $t=\f {1820}{9801}\sqrt {29}$ in Lemmas 3.2
and 3.3 and applying the above we deduce the result.

\pro{Theorem 3.9} Let $p$ be an odd prime such that $p\e
1,5,19,23\mod{24}$. Then
$$\sum_{k=0}^{p-1}\f{\b{2k}k^2\b{4k}{2k}}{28^{4k}}
\e\cases 4x^2\mod p&\t{if $p\e 1,19\mod{24}$ and so $p=x^2+18y^2$,}
\\0\mod{p^2}&\t{if $p\e 5,23\mod {24}$.}
\endcases$$
\endpro
Proof. From [LM, Table II] we know that the elliptic curve defined
by the equation $y^2=x^3+4x^2+(2-\f {40}{49}\sqrt{6})x$ has complex
multiplication by the order of discriminant $-72$. Thus, by (3.4)
 we have
$$\aligned&\sum_{n=0}^{p-1}\Ls{n^3+4n^2+(2-\f {40}{49}\sqrt{6})n}p
\\&=\cases 2x&\t{if $p\e 1,19\mod {24}$ and so $p=x^2+18y^2$,}
\\0&\t{if $p\e 5,23\mod {24}.$}\endcases\endaligned$$
Now taking $m=28^4$ and $t=\f {20}{49}\sqrt {6}$ in Lemmas 3.2 and
3.3 and applying the above we deduce the result.
\par\q
\newline{\bf Remark 3.4} Let $p>7$ be a prime. Z.W. Sun conjectured
that ([Su1, Conjecture A28])
$$\sum_{k=0}^{p-1}\f{\b{2k}k^2\b{4k}{2k}}{28^{4k}}
\e\cases 4x^2-2p\mod p&\t{if $p=x^2+2y^2\e 1,3\mod 8$,}
\\0\mod{p^2}&\t{if $p\e 5,7\mod 8$.}
\endcases$$

\pro{Theorem 3.10} Let $p$ be an odd prime such that $p\e \pm 1\mod
5$. Then
$$\sum_{k=0}^{p-1}\f{\b{2k}k^2\b{4k}{2k}}{(-2^{14}\cdot 3^4\cdot 5)^k}
\e\cases 4x^2\mod p&\t{if $p=x^2+25y^2$,}
\\0\mod{p^2}&\t{if $p\e 3\mod 4$.}
\endcases$$
\endpro
Proof. From [LM, Table II] we know that the elliptic curve defined
by the equation $y^2=x^3+4x^2+(2-\f {161}{180}\sqrt{5})x$ has
complex multiplication by the order of discriminant $-100$. Thus, by
(3.4)  we have
$$\aligned\sum_{n=0}^{p-1}\Ls{n^3+4n^2+(2-\f {161}{180}\sqrt{5})n}p
=\cases 2x&\t{if $p=x^2+25y^2$,}
\\0&\t{if $p\e 3\mod 4.$}\endcases\endaligned$$
Now taking $m=-2^{14}\cdot 3^4\cdot 5$ and $t=\f {161}{360}\sqrt{5}$
in Lemmas 3.2 and 3.3 and applying the above we deduce the result.
\par\q
\newline{\bf Remark 3.5} Let $p>7$ be a prime. Z.W. Sun made a conjecture
([Su1, Conjecture A25]) equivalent to
$$\sum_{k=0}^{p-1}\f{\b{2k}k^2\b{4k}{2k}}{(-2^{14}\cdot 3^4\cdot 5)^k}
\e\cases 4x^2-2p\mod p&\t{if $p=x^2+25y^2$,}
\\-4xy\mod{p^2}&\t{if $p=x^2+y^2$ with $5\mid x-y$,}
\\0\mod{p^2}&\t{if $p\e 3\mod 4$.}
\endcases$$

\pro{Theorem 3.11} Let $p>3$ be a prime. Then
$$\align&\sum_{k=0}^{p-1}\f{\b{2k}k^2\b{4k}{2k}}{648^k}
\e 0\mod{p^2}\qtq{for}p\e 3\mod 4,
\\&\sum_{k=0}^{p-1}\f{\b{2k}k^2\b{4k}{2k}}{(-144)^k}
\e 0\mod{p^2}\qtq{for}p\e 2\mod 3,
\\&\sum_{k=0}^{p-1}\f{\b{2k}k^2\b{4k}{2k}}{(-3969)^k}
\e 0\mod{p^2}\qtq{for}p\e 3,5,6\mod 7.\endalign$$
\endpro
Proof. This is immediate from Corollary 2.2 and [S2, Theorems
4.3-4.5].
\par \q
\par We remark that Theorem 3.11 was conjectured by the author in
[S1].

\enddocument
\bye